\documentclass[12pt,a4paper]{article}
\usepackage[english]{babel}
\usepackage{amsmath}
\usepackage{amstext}
\usepackage{amssymb}
\usepackage{amsthm}
\usepackage{amsfonts}
\usepackage{fullpage}
\usepackage{graphicx}
\usepackage{epsf}
\usepackage[dvips]{epsfig}
\usepackage{stmaryrd}
\usepackage[T1]{fontenc}
\usepackage{latexsym}
\usepackage{psfrag}

\title{Quasiperiodic infinite words : multi-scale case and dynamical properties}
\author{Thierry Monteil\footnote{Institut de Math\'ematiques de Luminy, CNRS UMR 6206,
Case 907, 163 Avenue de Luminy, 13288 Marseille cedex 09, France
-- monteil@iml.univ-mrs.fr}
\ \ and Solomon Marcus\footnote{Romanian Academy, Mathematics, Calea Victoriei 125, Bucharest, Romania 
-- Solomon.Marcus@imar.ro}}

\date{}

\newcommand{\N}{\mathbb N}
\newcommand{\R}{\mathbb R}

\newcommand{\defi}{\stackrel{\mbox {\tiny def}}{=}}

\newenvironment{appli}{\left( \begin{array}{ccc}}{\end{array} \right)}
\newcommand{\dans}{& \longrightarrow &}
\newcommand{\donne}{& \longmapsto &}
\newcommand{\ba}{\begin{appli}}
\newcommand{\ea}{\end{appli}}

\newenvironment{demo}{\noindent {\bf Proof:}}{\hfill $\square$ \\}
\newtheorem{thm}{Theorem}[]
\newtheorem{cor}{Corollary}[]

\newcommand{\tq}{\ | \ }

\renewcommand{\[}{[ \, }
\renewcommand{\]}{ \, ]}

\newcommand{\khi}{\chi}
\newcommand{\card}{\mbox{card}}

\newcommand{\cge}{\xrightarrow[ n \rightarrow \infty]{}}

\begin{document}

\maketitle


\vspace{2 cm}

\begin{abstract}

\noindent An infinite word $x$ is said to be quasiperiodic if there exists a finite
word $q$ such that $x$ is covered by occurrences of $q$ (such a $q$ is
called a quasiperiod of $x$). Using the notion of derivation, we show that this definition is not
sufficient to imply any symmetry in an infinite word. 
Therefore we introduce multi-scale quasiperiodic words, i.e. quasiperiodic
words that admit an infinite number of quasiperiods. Such words are
uniformly recurrent, this allows us to study the subshift they generate.
We prove that multi-scale quasiperiodic subshifts
are uniquely ergodic and have zero topological entropy as well as zero Kolmogorov complexity. 
Sturmian subshifts are shown to be multi-scale quasiperiodic.\\ 



\end{abstract}

\clearpage

\section{Introduction and notations}

\subsection{Background: quasiperiodicity and symmetry}

The general concept of symmetry 
is central in the study of infinite
words, and many notions are devoted to measure it, like word complexity, 
recurrence, Kolmogorov's complexity, entropy, {\ldots}.
Those notions are more or less linked and in any case the most symmetric
words are the periodic ones.

In \cite{Mar1}, the second author introduced a new notion to describe symmetric words: quasiperiodicity. 
This notion is an extension of the notion of quasiperiodicity for finite 
words that was introduced in \cite{ApoEhr} 
and has also roots in different contexts like musicology \cite{CrIlRa} and molecular biology \cite{KMGL} \cite{MilJur}.\\

Let $A$ be a finite set (called {\em alphabet}).
Let $x=x_0 x_1 x_2 \dots$ in $A^{\N}$ be an infinite word. We will say that $x$ is {\em
quasiperiodic} if there exists a finite word $q$ in $A^*$ such that $x$ is
covered by the occurences of $q$ ($A^*$ denotes the set of finite words on the alphabet $A$).
Such a $q$ is called a {\em quasiperiod} of $x$ ; it is a prefix of $x$.

More precisely, let $(i_n)_{n\in \N}$ be the increasing sequence whose image is the set 
$\{i\in\N \tq x_{i\rightarrow i+l(q)-1}=q\}$, where $l(q)$ denotes the length of $q$, 
and $x_{i\rightarrow j}$ denotes $x_i x_{i+1} x_{i+2} \dots x_j$
(for $0\leq i\leq j$).
The word $q$ is a quasiperiod of $x$ whenever $(i_n)$ is a well defined infinite sequence such that 
$i_0=0$ and $i_{n+1}-i_{n}\leq l(q)$ for any integer $n$.

For example, the infinite word $x=ababaabaabaababababaabababaabaabaababaaba{\ldots}$ 
is quasiperiodic with $aba$ as a quasiperiod 
and the infinite quasiperiodic words that admit $aba$ as a quasiperiod are precisely 
the words over the alphabet $\{a,b\}$ beginning with the word $ab$ 
and in which the words $aaa$ and $bb$ do not appear.\\

The second author was then looking for relations with other notions of symmetry. 
In \cite{Mar1} and \cite{LevRic1}, it is shown that there is no such relation. 
Let us introduce a basic tool that will unify those results and 
allow us to understand this lack of fitting between quasiperiodicity and other notions: the derivation.

\subsection{The derivation: a change of scale}
\subsubsection{Derivation}

Let $x=x_0 x_1 x_2 \dots$ be a quasiperiodic word, and $q$ be a quasiperiod of $x$.
Let $(i_n)_{n\in \N}$ be the sequence as defined before.
We can define $\frac{\partial x}{\partial q}$ to be the infinite word on the alphabet $\{0,\dots,l(q)-1\}$
whose $k^{\mbox{\tiny th}}$ letter is $l(q)-i_{k+1}+i_k$ ($k\geq 0$).

In other terms, the $k^{\mbox{\tiny th}}$ letter of $\frac{\partial x}{\partial q}$ is the length of the overlap
between the $k^{\mbox{\tiny th}}$ and the $(k+1)^{\mbox{\tiny th}}$ occurrence of $q$ in $x$.\\

For example, if $x=ababaabaabaababababaabababaabaabaababaaba{\ldots}$ is quasiperiodic with $aba$ as a quasiperiod, 
then $\frac{\partial x}{\partial aba} = 100011101100010{\ldots}$ .

\subsubsection{Integration}
In some sense, the $\frac{\partial}{\partial q}$ operator removes exactly the information contained in $q$ 
since the knowledge of $\frac{\partial x}{\partial q}$ and $q$ is sufficient to reconstruct $x$.

Indeed, we can consider the reverse operation : the integration. 
If $x$ is an infinite word on a finite alphabet $A\subset \N$ and if $w$ is a finite
word whose length is greater than $\max A$, 
then we can define the word $\int_w x$ as the image of $x$ under the substitution $\sigma_w$ 
that replaces the occurences of $i\in A$ by the $l(w)-i$ first letters of $w$.\\

For example, if $x=01121010201...$ and $w=aabcaa$, 

then  $\int_w x=aabcaaaabcaaabcaaabcaabcaaabcaaaabcaaabcaaaabcaabcaaaabcaaabcaa...$.\\

If $w$ is such that for each $n$ in $A$, $w$ is a prefix of $\sigma_w(n).w$, 
then $\int_w x$ is quasiperiodic with $w$ as a quasiperiod, no matter how $x$ is random.
This is the case for example if $A\subset \{0,\dots,n\}$ and $w=a^nba^n$ ($n\geq 0$).

\subsubsection{Quasiperiodicity is spread everywhere}

\begin{description}

\item[Recurrence] 
An infinite word $x$ is said to be {\em recurrent} if any finite word $u$ appearing in $x$ appears infinitely often in $x$.
It is said to be {\em uniformly recurrent} if moreover for each finite word $u$ appearing in $x$,
the gap between two consecutive occurences of $u$ in $x$ is bounded.

To construct a non recurrent quasiperiodic word, just take a non recurrent word $x$ on the alphabet $\{0,1\}$
and consider $\int_{aba} x$.
Do the same to construct a uniformly recurrent quasiperiodic word. 

Minimality is the analogue of uniform recurrence in the vocabulary of topological dynamical systems 
(see section \ref{section-minimal}).

\item[Complexity]
For any infinite word $x$ and any integer $n$, $L_n(x)$ denotes the finite words of length $n$ that occur in $x$ 
and we define $L(x)\defi \bigcup_{n\in \N} L_n(x)$.

The function that sends an integer $n$ to $p_n(x) \defi \card(L_n(x))$ is called the {\em word complexity} of $x$.

We can define an equivalence relation on the set of infinite words 
through the asymptotic behaviour of their complexity function: 
two words $x$ and $y$ are said to be {\em complexity equivalent} if there exists a positive integer K such that
for all $n\geq 1$, $p_n(x)\leq K p_{Kn}(y)$ and $p_n(y)\leq K p_{Kn}(x)$.
Hence, bounded, linear, quadratic, polynomial or exponential growths are preserved under complexity equivalence.

\begin{thm}
There are quasiperiodic words in any class of complexity equivalence.
\end{thm}

\begin{demo}
Let $x$ be an infinite word on an alphabet $A$. There is no restriction to suppose that $A=\{0,\dots,k\}$.
Let $y\defi \int_{a^kba^k} x$. We check that $y$ is a quasiperiodic word in the class of complexity equivalence of $x$.
\end{demo}

The non ultimately periodic words with the smallest word complexity are the sturmian words: 
an infinite word is said to be {\em sturmian} if for any integer $n$, $p_n(x)=n+1$.
In \cite{LevRic1}, Florence Lev\'e and Gw\'ena\"el Richomme proved that there exists sturmian
words that are not quasiperiodic.

Entropy is the exponent of the complexity function in the vocabulary of topological dynamical systems 
(see section \ref{section-complexity}).

\item[Frequencies]
If $u$ and $v$ are finite words, let $\#(u,v)$ denotes the number of occurences of $u$ in $v$.
An infinite word $x$ is said to have {\em frequencies} if for any finite word $w$,
$\frac{1}{n}\#(w, x_{0\rightarrow n-1})$ admits a limit when $n$ tends to infinity.

By integration, there exists quasiperiodic words that do not have frequencies: 
if $x$ is an infinite word such that $0$ does not appear with frequencies, 
then $\int_{aba} x$ is a quasiperiodic word such that $bab$ does not appear with frequencies.

Unique ergodicity is a strong analogue of having frequencies in the vocabulary of topological dynamical systems 
(see section \ref{section-ue}).\\
 
\end{description}

Hence the notion of quasiperiocicity does not insert well among other notions of symmetry.

\subsection{Multi-scale quasiperiodic words}

We have to notice that all classical notions of symmetry are invariant under such a renormalization procedure 
(derivation corresponds to induction if we are studying dynamical properties like entropy \cite{Abr}).
Therefore, if we want that a notion says something about symmetry, we should ensure that it is stable under such a change of scale.
As we saw with the derivation procedure, 
the existence of a quasiperiod $q$ in a word $x$ just imposes rigidity at the scale around $l(q)$ 
but does not impose anything at larger scales.

This leads to the following definition : 
an infinite word is said to be {\em multi-scale quasiperiodic} if the set $Q(x)$ of its quasiperiods is infinite.\\

The easiest non-periodic multi-scale quasiperiodic words we can construct, are the fixed points for some particular integration operators. 
For example, the fixed point of $\int_{010}$ is multi-scale quasiperiodic, it is known as the Fibonacci word 
(a precise description of the quasiperiods of this word can be found in \cite{LevRic1}). 
There exists much wilder multi-scale quasiperiodic words (see Theorem \ref{thm-sup-complexity} in section \ref{section-complexity}).\\

In the next section we will prove that multi-scale quasiperiodic words are uniformly recurrent. 
This will allow us to study the subshift generated by them (section \ref{section-minimal}). 
Concerning the complexity, we will prove that multi-scale quasiperiodic subshifts have zero topological entropy 
as well as zero Kolmogorov complexity (section \ref{section-complexity}).
We will also prove that sturmian subshifts are multi-scale quasiperiodic (section \ref{section-sturmian}). 
Concerning frequencies, we will prove that multi-scale quasiperiodic subshifts are uniquely ergodic (section \ref{section-ue}).

\clearpage

\section{Uniform recurrence and minimality}\label{section-minimal}

An infinite word $x\in A^{\N}$ is said to be {\em uniformly recurrent} if 
any finite word $u\in L(x)$ occurs infinitely many times in $x$
and the gap between two consecutive occurrences of $u$ in $x$ is bounded, 
equivalently

$$\forall u \in L(x) \ \ \ \exists n \geq 1 \ \ \ \forall v \in L_n(x) \ \ \ \#(u,v)\geq 1$$

\begin{thm}\label{thm-unif-rec}
Any multi-scale quasiperiodic word $x$ is uniformly recurrent.
\end{thm}

\begin{demo}
Let $u$ be a finite word that occurs in $x$.
Since every quasiperiod of $x$ is a prefix of $x$, one of them must
contain an occurrence of $u$ (they have unbounded length). Let $q$ be such a
quasiperiod. Any word in $L_{2l(q)}$ contains at least an occurrence of $q$
and therefore at least an occurrence of $u$.
\end{demo}

This property of multi-scale quasiperiodic words allows us to deal with the
dynamical system generated by a multi-scale quasiperiodic word as follows:\\

We endow $A$ with the discrete topology and $A^\N$ with the product
topology. This makes $A^\N$ a metrisable compact space.

We note $$ S \defi \ba A^\N \dans A^\N \\ 
x=x_{0} x_{1} \dots x_{n} \dots \donne 
x_{1} x_{2} \dots x_{n+1} \dots \ea$$ for the shift. It is a continuous
map.

If $x$ is a multi-scale quasiperiodic word, we define 
$$X \defi \overline{\{ S^k(x) \tq k\in \N \}}$$ 
and we still note $S$ for the restriction of $S$ to $X$, making $(X,S)$ a
topological dynamical system.

Theorem \ref{thm-unif-rec} is equivalent to say that $(X,S)$ is a {\em
minimal subshift} i.e. $X$ is a nonempty closed subset of $A^\N$ 
stable under $S$ and that is minimal for those properties.

A minimal subshift generated by a multi-scale quasiperiodic word is called a
{\em multi-scale quasiperiodic subshift}.

If $u$ is a finite word, we define the {\em cylinder}  

$$[u] \defi \{ x \in X \tq  (\forall i\leq n-1) (x_i=u_i) \} \ \ \ 
(u= u_0 \dots u_{n-1} \in A^n \ \ (n\in \N))$$

We can notice that the derivation of a multisacled quasiperiodic word over the quasiperiod $q$ 
corresponds to the induction of the subshift $(X,S)$ on the cylinder $[q]$.

We extend the notion of language and complexity to minimal subshifts:
$$L_n(X) \defi \{u \in A^n \tq [u]\neq \emptyset\} \ \mbox{ and } \ p_n(X)=\card
(L_n(X)) \ \ \ (n \in \N)$$

$$ L(X) \defi
\bigcup_{n\in \N} L_n(X)$$ 
If $y$ is in $X$, we have $L(y)=L(X)$ and $p_{.}(y)=p_{.}(X)$.

\clearpage

\section{Complexity and topological entropy}\label{section-complexity}

\subsection{Word complexity}

\begin{thm}\label{thm-inf-complexity}
Let $x$ be a multi-scale quasiperiodic word. 
Then $$\liminf_{n\rightarrow \infty} \frac{p_n(x)}{n^2} \leq 1 < \infty$$
\end{thm}

\begin{demo}
Let $q\in Q(x)$ and $u\in L_{l(q)}(x)$.
Since $q$ is a quasiperiod of $x$, $u$ is a subword of some $v.q$ where $v$
is nonempty prefix of $q$. $u$ is determined by the choice of
$v$ and his position in $v.q$. 
There are $l(q)$ prefixes of $q$ and for such a prefix $v$, there are $l(v)
\leq l(q)$ available positions for $u$ (we do not count $q$ several times).
Finally, there are only $l(q)^2$ possibilities for $u$ and $p_{l(q)}(x) \leq
l(q)^2$. The result follows since $\{ l(q) \tq q\in Q(x)\}$ is not bounded.
\end{demo}

\begin{thm}\label{thm-sup-complexity}
For each positive function $f:\N\rightarrow \R_+^*$ 
that converges to zero, 
there exists a multi-scale quasiperiodic subshift $(X,S)$
such that 
$\frac{1}{n} \log (p_n(X)) \geq f(n)$ for infinitely many $n$.
In particular, we can ask to $p_n(X)$ 
to grow faster than any polynomial on a subsequence.
\end{thm}

\begin{demo}
Let $\varphi:\N\rightarrow\N^*$ be an increasing sequence such that 
$2nf(2n\varphi(n))\leq 1$ for any $n$.
Let $A=\{0,1\}$ be the alphabet.
For $n\geq 1$, let $w_n$ be a finite word on $A$ such that 
\begin{itemize}
  \item $w_n$ begins and ends with the letter $0$,
  \item every word $w$ of length $2n$ such that $\#(0,w)=\#(1,w)=n$ appears as a subword of $w_n$.
\end{itemize}
Now, let us define by induction a sequence of finite words over $A$:
\begin{itemize}
  \item $u_0=010$
  \item $u_{n+1}=\int_{u_n}w_{\varphi(l(u_n))}$
\end{itemize}
which is well defined since for any $n$, $u_n$ begins and ends with the letter $0$.

Since $u_n=\int_{u_n}0$, $u_n$ is a prefix of $u_{n+1}$: 
let $x$ be the unique infinite word over $A$ such that 
$u_n$ is a prefix of $x$ for any $n$.

If $k\leq l$, $u_l$ is covered by occurences of $u_k$,
so all the $u_n$ are quasiperiods of $x$ and
$x$ is multi-scale quasiperiodic. 
Let $(X,S)$ denote the associated minimal subshift

Let $n$ be a positive integer. 
$u_{n+1}$ and therefore $x$ contain an occurence of $\int_{u_n}w$,
where $w$ is any word of length $2\varphi(l(u_n))$ such that 
$\#(0,w)=\#(1,w)=\varphi(l(u_n))$.
There are at least $2^{\varphi(l(u_n))}$ such different words 
and each of them has length $\varphi(l(u_n))(l(u_n)+(l(u_n)-1))$.

Hence, $$p_{2\varphi(l(u_n))l(u_n)}(X) 
\geq p_{\varphi(l(u_n))(l(u_n)+(l(u_n)-1))}(X)
\geq 2^{\varphi(l(u_n))}$$

So, $$\frac{1}{2\varphi(l(u_n))l(u_n)} \log(p_{2\varphi(l(u_n))l(u_n)}(X)) 
\geq \frac{\varphi(l(u_n))}{2\varphi(l(u_n))l(u_n)} = \frac{1}{2l(u_n)}
\geq f(2\varphi(l(u_n))l(u_n))$$
\end{demo}

\clearpage

We can notice that for the infinite words constructed here, 
the scales controlled by the quasiperiods are very sparse.

\subsection{Topological entropy}

If $(X,S)$ is a minimal subshift, then the limit 
$$h_{top}(X) \defi \lim_{n\rightarrow \infty} \frac{1}{n} \log (p_n(X))$$ 
exists and is named the {\em topological entropy} of $X$.

\begin{cor}
Any multi-scale quasiperiodic subshift has zero topological entropy.
\end{cor}

\subsection{Kolmogorov complexity}

Let $U$ be a fixed universal Turing machine and for each finite word $u\in A^*$,
let $K_U(u)$ denotes the Kolmogorov complexity associated to $u$, 
i.e. the length of the shortest binary word $p$ such that $U(p)=u$
(see \cite{Bru}).
For a minimal subshift $(X,S)$, we can define 
$$K(X) \defi \sup_{x\in X} \limsup_{n\rightarrow \infty} \frac{K_U(x_{1\rightarrow n})}{n}$$
(this number is independant of the choice of $U$ since if $U'$ is another universal Turing machine, 
there is a constant $C$ such that for any finite word $u$, $K_U(u)-C  \leq K_{U'}(u)  \leq  K_U(u)+C$).

\begin{cor}
For any multi-scale quasiperiodic subshift $(X,S)$, $K(X)=0$.
\end{cor}

Indeed, \cite{Bru} Theorem 3.1 asserts that $K(X)\leq h_{top}(X)$. \\

But we can also give a more direct argument : 
if $q$ is a fixed quasiperiod of a multi-scale quasiperiodic subshift $(X,S)$ and if $x$ is in $X$, 
then $\limsup_{n\rightarrow \infty} \frac{K_U(x_{1\rightarrow n})}{n} \leq \frac{4\log(l(q))}{l(q)}$.
To prove this, it suffices to remark that the integration algorithm can be coded in $O(1)$.

Then, if $n$ is bigger than $4l(q)$, there exists three finite words $u$, $w$ and $v$ 
such that $x_{1\rightarrow n}=u.w.v$, the lengths of $u$ and $v$ are uniformly bounded by $l(q)$,
and $w$ begins and ends with $q$ (and can therefore be derivated).

The remaining problem is the control of the length of $\partial w / \partial q$:
if the occurences of $q$ overlap each other deeply, $\partial w / \partial q$ can be rather long.
To solve this, we decide to replace recursively an occurence of $n_1.n_2$ in $\partial w / \partial q$ by $n_1 + n_2$ 
if $n_1$ and $n_2$ are smaller than $n/2$.
This operation 
consists in omitting
some useless occurences of $q$ in $x_{1\rightarrow n}$.
Then at least half of the numbers appearing in the new form of $\partial w / \partial q$ are bigger than $l(q)/2$, 
so the length of this new $\partial w / \partial q$ will be less than $4n/l(q)$.

Since the coding of each letter of $\partial w / \partial q$ costs $\log(l(q))$, we have a total cost less than 
a constant (to code the integration algorithm, $q$, $u$ and $v$) plus $\log(l(q)) 4n/l(q)$ 
(to code the new form of $\partial w / \partial q$).



\clearpage

\section{Unique ergodicity and frequencies}\label{section-ue}

For a multi-scale quasiperiodic subshift $(X,S)$, 
we will now study the set $\mathcal{M}(X,S)$ of Borel
probability measures on $X$ that are invariant under $S$.
This set can be identified with a nonempty compact convex subset of $C^0(X,\R)'$
endowed with the weak-star topology. 

A $S$-invariant measure $\mu \in \mathcal{M}(X,S)$ is said to be {\em ergodic} 
if the only Borel sets $A \subset X$ such that $S^{-1}(A)=A$ have measure
$\mu(A) = 0 \ \mbox{ or } \ 1$. Such measures are the extremal points of
$\mathcal{M}(X,S)$ and satisfies Birkhoff's theorem :

$$\forall f\in L^1(X,\R) \ \ \ \ \frac{1}{n} \sum_{k=0}^{n-1}
f \circ S^k \xrightarrow[ n \rightarrow \infty]{\mu - a.e.}
\int_{X} f d\mu $$

A minimal subshift is said to be {\em uniquely ergodic} if 
$\card( \mathcal{M}(X,S) ) = 1$.

One interest of such a situation is that, the unique invariant measure
$\mu$ is ergodic, moreover the convergence in Birkhoff's theorem is uniform
for continuous functions.

\begin{thm}\label{thm-UE}
Any multi-scale quasiperiodic subshift $(X,S)$ is uniquely ergodic.
\end{thm}

\begin{demo}

We will first construct a $S$-invariant probability measure on $X$ and then
prove that it is the only one. 

{\bf Step 1:} We construct a candidate to be the unique measure.
For this, we will approximate $X$ by periodic
subshifts generated by the $q^{\omega}=qqq\dots$ for $q\in Q(X)$.  

For $q$ in $Q(X)$, let $$\mu_q \defi \frac{1}{l(q)} \sum_{k=0}^{l(q)-1} \delta_{S^k(q^{\omega})}$$

($\delta$ stands for the one-point Dirac's measure).

$\mu_q$ is the only element of $\mathcal{M}(A^{\N},S)$ that gives
mesure $1$ to the periodic subshift generated by the periodic word $q^{\omega}$.
By compacity we can find an infinite subset $Q' \subset Q(X)$ such that 
$$\mu_q \xrightarrow[ \stackrel{l(q) \rightarrow \infty}{q\in Q'}]{} \mu$$
for some $\mu$ in $\mathcal{M}(A^{\N},S)$.

Note that if $X$ is aperiodic the $\mu_q$'s give measure $0$ to $X$. However $\mu$ will give strictly 
positive measure to $X$ (as we will see in Step 2) : it shouldn't be surprising 
since the characteristic function of $X$ is not
continuous.\\ 

{\bf Step 2:} Let us show that $\mu(X)=1$.
Since $X$ is closed, we have the following approximation by clopen
sets: $$X = \overline{X} = \bigcap_{n\geq 1} \bigcup_{u\in L_n(X)} [u].$$

Let $n\geq 1$ and let $q\in Q'$ such that $l(q)\geq n$. For $i\in \{ 0,
\dots, l(q)-n \}$, we have $q^{\omega}_{i \rightarrow i+n-1} \in L_n(X)$ (as a subword of
$q$).
Hence $\mu_q(\bigcup_{u\in L_n(X)} [u]) \geq (l(q)-n+1)/l(q)$.

Letting $l(q)$ tending to infinity, since the caracteristic function
of $\bigcup_{u\in L_n(X)} [u]$ is continuous, we have
$\mu(\bigcup_{u\in L_n(X)} [u]) = 1$. By countable intersection
($n$ is arbitrary), we have $\mu(X)=1$.
Hence, we can still denote by $\mu$ for the restriction of $\mu$ to
$X$. \\ 

{\bf Step 3:} Let $\nu$ be an ergodic measure on $X$. We will show
that $\nu = \mu$

By Birkhoff's theorem, there is $x$ in $X$ such that for $u\in
L(X)$, 
$$\nu([u])=\lim_{n\rightarrow \infty} \frac{1}{n} \#(u,x_{0 \rightarrow n+l(u)-2})=
\lim_{n\rightarrow \infty} \frac{1}{n} \#(u,x_{0 \rightarrow n-1})$$

Let $q$ in $Q'$ such that $l(q)\geq l(u)$. We decompose $x$ into blocks of length $2l(q)$:
$x=B_0.B_1.B_2.B_3.B_4{\dots}$ with $B_i=x_{2l(q)i \rightarrow 2l(q)(i+1)-1}$. 
Since each $B_i$ is in $L_{2l(q)}(X)$, it contains at least one occurence
of $q$, hence $\#(u,B_0.B_1.B_2.{\cdots}.B_i)\geq (i+1) \#(u,q)$ for each
$i$ in $\N^*$.

We have 
$$\nu([u])=\lim_{i\rightarrow \infty}\frac{1}{2l(q)i} \#(u,B_0.B_1.B_2.{\cdots}.B_i) 
\geq \lim_{i\rightarrow \infty}\frac{i+1}{2l(q)(i+1)} \#(u,q)
= \frac{1}{2l(q)}\#(u,q) $$

Moreover, we can control the frequency of occurrences of $u$ in
$q^{\omega}$ by estimating the occurrences of $u$ in $q$ and 
bounding the number of 
occurrences of $u$ that appear between two consecutive
occurences of $q$ in $q^{\omega}$:
$$\mu_q([u])=\lim_{n\rightarrow \infty} \frac{1}{n} \#(u,q^w_{0 \rightarrow
n+l(u)-2})\leq \frac{1}{l(q)}(\#(u,q)+l(u)) = \frac{1}{l(q)}\#(u,q) +\frac{l(u)}{l(q)} $$

Therefore, $$\mu_q([u]) \leq \frac{1}{l(q)}\#(u,q) + \frac{l(u)}{l(q)}  
\leq 2\nu([u]) + \frac{l(u)}{l(q)} $$

Letting $l(q)$ tend to infinity, we have $\mu([u]) \leq 2 \nu([u])$. So,
$\mu$ is absolutely continuous relatively to $\nu$. 
It is well known that this implies $\mu=\nu$, but for sake of completeness, 
we include a short proof here.
There exists a
measurable function $f\in L^1(X,\R_+)$ such that for each borel set $A\subset X$, 
$\mu(A)= \int_{A}^{} f d\nu$. Since $\mu$ is $S$-invariant, we have
$\int_{A}^{} f d\mu=\int_{S^{-1}A}^{} f d\mu$ for each borel set $A\subset X$.

Let us show that $f$ is constant almost everywhere. Assume by contradiction
that the measure of set $A \defi \{ x \in X / f(x) \geq \int_{X}^{} f d\nu
\}$ is in $\]0,1\[$. Since $\nu$ is $S$-ergodic, $A$ is not $S$-invariant, so
$\nu (S^{-1}(A) \setminus A) =  \nu (A \setminus S^{-1}(A))  > 0$.

Hence 
$ \nu(A\setminus S^{-1}(A)) \int_{X}^{} f d\nu 
\leq \int_{A \setminus S^{-1}(A)} f d\nu 
= \int_{A} f d\nu -  \int_{A \cap S^{-1}(A)} f d\nu =  
\int_{S^{-1}(A)} f d\nu -  \int_{A \cap S^{-1}(A)} f d\nu =
\int_{S^{-1}(A)\setminus A} f d\nu 
< \nu(A \setminus S^{-1}(A)) \int_{X}^{} f d\nu $ which is absurd.

Hence $f$ is constant with value $\nu(X)=1$, so $\mu=\nu$ and $\mu$ is the
only $S$-invariant measure.

\end{demo}

\begin{cor}
Let $x$ be a multi-scale quasiperiodic word. Then each finite word $u$
occuring in $x$ has frequencies i.e. $\frac{1}{n} \#(u,x_{0\rightarrow n-1})$
converges when $n\rightarrow \infty$.
\end{cor}

\begin{demo}
Let $\mu$ be the unique $S$-invariant measure for the associated subshift
$(X,S)$. The caracteristic function $\khi_{[u]}$ of $[u]$ is continuous, so we have a
uniform and therefore a pointwise convergence in Birkhoff's theorem : 
$$ \frac{1}{n} \#(u,x_{0\rightarrow n-1}) =  \frac{1}{n} \sum_{k=0}^{n-1}
\khi_{[u]}(S^k(x)) \cge \mu([u])$$
\end{demo}

\clearpage

We can remark that this proof has the same flavour as the the main result of
\cite{Bos} that asserts that every subshift with subaffine complexity has only
a finite number of ergodic measures. In fact, Theorem \ref{thm-UE} and
Boshernitzan's result can both be deduced from a more general statement that
involves the geometry of Rauzy graphs associated to a minimal subshift:\\

To each minimal subshift $(X,S)$ we can associate a sequence $(G_n)_{n\geq 1}$ of
oriented graphs as follows: the vertices of $G_n$ is $L_n(X)$ and there is an
edge from $u$ to $v$ if and only if there exists $w$ in $L_{n+1}(X)$ such
that $u$ is a prefix of $w$ and $v$ is a suffix of $w$.

Those graphs are named the {\em Rauzy graphs} associated to $(X,S)$.\\

If $K\geq 1$, $(X,S)$ is said to be $K${\em -deconnectable} if there exists an
extraction $\alpha \in \ \uparrow \!\! (\N^*,\N^*)$ and a constant $K'\geq1$ 
such that for all $n\geq 1$ 
there exists a subset $D_{\alpha(n)} \subset L_{\alpha(n)}(X)$ of at most $K$ vertices 
such that every path in $G_{\alpha(n)}(X) \setminus D_{\alpha(n)}$ 
is of length at most $K' \alpha(n)$ (in particular it do not contains any cycle).
This means that, up to extraction, we can disconnect (in a specific way) the Rauzy graphs 
by removing at most $K$ vertices.

\begin{thm}[\cite{Mon}]
A $K$-deconnectable minimal subshift has at most $K$ $S$-invariant ergodic measures.
\end{thm}

This result implies Boshernitzan's one, by taking for $D_n$ the set of right
special factors (i.e. the set of vertices having outgoing degree strictly greater
than one). It also implies Theorem \ref{thm-UE} since every multi-scale quasiperiodic
subshift is $1$-deconnectable. Indeed, if $q$ is in $Q(X)$,
$G_{l(q)}(X) \setminus \{q\}$ does not contain any path of length
greater than $l(q)$.

%
%
%

\clearpage

\section{Sturmian subshifts are multi-scale quasiperiodic}\label{section-sturmian}
In \cite{LevRic1}, Florence Lev\'e and Gw\'ena\"el Richomme proved that there exists sturmian
words that are not quasiperiodic. 
In terms of complexity, sturmian words are the more symmetric words after periodic ones, 
so this result do not comfirm that quasiperiodicity fits well with other notions of symmetry.

The dynamical point of view will allow us to solve the problem:

\begin{thm}
Sturmian subshifts are multi-scale quasiperiodic.
\end{thm}

\begin{demo}
Let $(X,S)$ be a sturmian subshift.
Since $p_1(X)=2$, we can consider that $X$ is defined on the alphabet $\{a,b\}$.

A word $u\in L_n{X}$ is said to be {\em left special} if $au$ and $bu$ are in $L_{n+1}(X)$.
Since $p_{n+1}(X)-p_{n}(X)=1$, there exists exactly one left special word $l_n$ of length $n$ ($n\geq 0$).
A prefix of a left special word is still a left special word, so $l_{n+1}$ begins with $l_n$: 
let us denotes by $x$ the infinite word that begins by $l_n$ for any integer $n$. 
Since $L(x)\subset L(X)$, $x$ is in $X$ (remember that $X$ is closed). 
We will prove that $x$ is multi-scale quasiperiodic.\\

The evolution of the Rauzy graphs of sturmian subshifts is described by Rauzy (see \cite{ArnRau}): 
for infinitely many $n$ (named ``bursts''), the Rauzy graph $G_n(X)$ is eight shaped i.e. 
$G_n(X)$ is the union of two disjoints loops from $l_n$ to $l_n$.
By minimality, the minimal size of the two loops tends to infinity with $n$, 
in particular, it is positive for $n$ big enough.
Since the sum of the lengths of the two loops is equal to $p_n(X)=n+1$, 
then each loop has size less or equal than $n$.\\

So, for infinitely many $n$, any path in $G_n(X)$ of length $n$ starting from $l_n$ has to meet $l_n$ again i.e.
$l_n$ is a quasiperiod of $x$ ($x$ can be viewed as an infinite path in $G_n(X)$ starting from $l_n$).
Therefore, $(X,S)$ is a multi-scale quasiperiodic subshift.

\end{demo}

Recently, Florence Lev\'e and Gw\'ena\"el Richomme gave a precise description of the quasiperiodic sturmian words (\cite{LevRic2}).


\clearpage

\section{Conclusion}
The gain of fitting with other symmetry classes obtained by considering multi-scale quasiperiodic subshifts 
instead of quasiperiodic words can be sumarized as follow : 


\begin{figure}[h!]
\begin{center}
\psfrag{compl}{$\displaystyle \frac{\log(p_n(.))}{n} \rightarrow 0$}
\psfrag{unif}{uniformly recurrent}
\psfrag{quasip}{{quasiperiodic}}
\psfrag{periodic}{periodic}
\psfrag{sturm}{sturmian}
\psfrag{freq}{frequencies}
\psfrag{random}{random}
\includegraphics[width=.7 \linewidth]{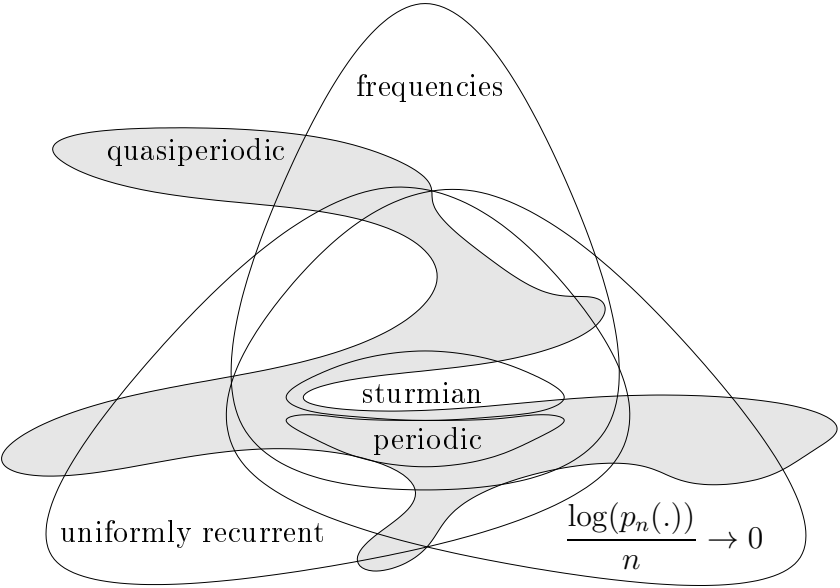}
\caption{\label{Bilan-quasiperiodic} Quasiperiodic words among other words.}
\end{center}
\end{figure}

\begin{figure}[h!]
\begin{center}
\psfrag{htop}{$h_{top}(.)=0$}
\psfrag{minimal}{minimal}
\psfrag{multiscaled}{multi-scale}
\psfrag{quasip}{{quasiperiodic}}
\psfrag{periodic}{periodic}
\psfrag{sturm}{sturmian}
\psfrag{ue}{uniquely ergodic}
\psfrag{random}{random}
\includegraphics[width=.7 \linewidth]{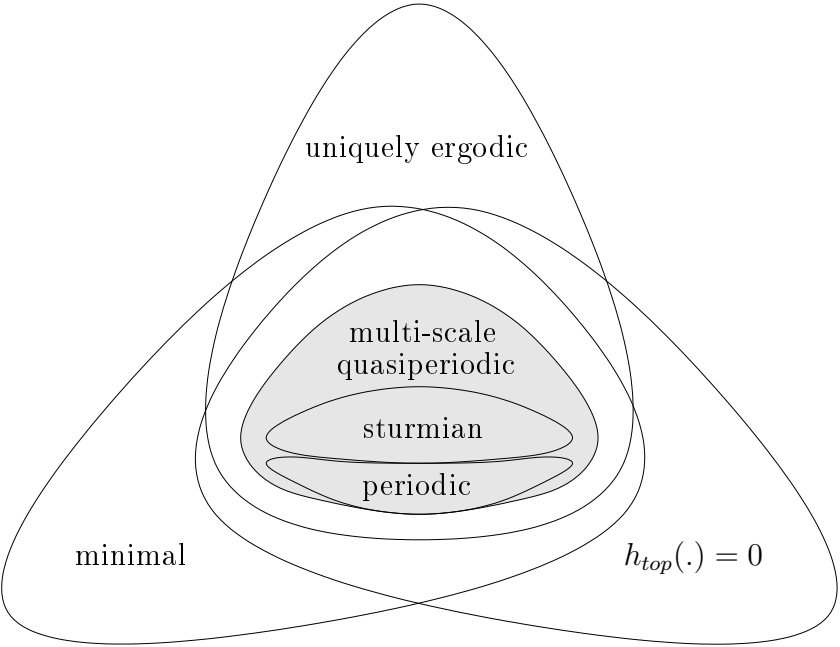}
\caption{\label{Bilan-multi-scale} multi-scale quasiperiodic subshifts among other subshifts.}
\end{center}
\end{figure}

%
%
%
%
%
%


\clearpage

\begin{thebibliography}{99}

\bibitem[Abr]{Abr}  L.M. {\sc Abramov},
{\it The entropy of an induced transformation}, Doklady Akad. Nauk SSSR 128 (1959), p 647--650


\bibitem[ApoEhr]{ApoEhr} A. {\sc Apostolico}, A. {\sc Ehrenfeucht},
{\it Efficient detection of quasiperiodicities in strings},
Theoret. Comput. Sci. 119 (1993), no. 2, p 247--265.


\bibitem[ArnRau]{ArnRau} P. {\sc Arnoux}, G. {\sc Rauzy},
{\it Représentation géométrique de suites de complexité $2n+1$}, 
Bull. Soc. Math. France 119 (1991), no. 2, p 199--215.


\bibitem[Bos]{Bos}  M. {\sc Boshernitzan},
{\it A unique ergodicity of minimal symbolic flows with linear block growth},
J. Analyse Math. 44 (1984/85), p 77--96.

\bibitem[Bru]{Bru} A.A. {\sc Brudno},
{\it Entropy and the complexity of the trajectories of a dynamic system},
Trudy Moskov. Mat. Obshch. 44 (1982), p 124--149.

\bibitem[CrIlRa]{CrIlRa} T. {\sc Crawford}, C. S. {\sc Iliopoulos}, R. {\sc Raman},
{\it String Matching Techniques for Musical Similarity and Melodic Recognition},
Computing in Musicology, vol. 11, p 73--100, 1998.

\bibitem[KMGL]{KMGL} S. {\sc Karlin}, M. {\sc Morris}, G. {\sc Ghandour}, M. {\sc Leung},
{\it Efficient algorithms for molecular sequence analysis},
Proc. National. Acad. Sci U.S.A. 85 (1988), p 841--845.


\bibitem[LevRic1]{LevRic1} F. {\sc Lev\'e}, G. {\sc Richomme},
{\it Quasiperiodic infinite words: some answers},
Bull. Eur. Assoc. Theor. Comput. Sci. EATCS No. 84, (2004), p 128--138.


\bibitem[LevRic2]{LevRic2} F. {\sc Lev\'e}, G. {\sc Richomme},
{\it Private communication, april 2005}.

\bibitem[Mar1]{Mar1}  S. {\sc Marcus},
{\it Quasiperiodic infinite words},
Bull. Eur. Assoc. Theor. Comput. Sci. EATCS No. 82, (2004), p 170--174.


\bibitem[Mar2]{Mar2}  S. {\sc Marcus},
{\it Symmetry phenomena in infinite words, with biological, philosophical and aesthetic relevance},
\texttt{http://symmetry.hu/contents/contents-2003-2004.pdf}


\bibitem[MilJur]{MilJur} A. {\sc Milosavljevic}, J. {\sc Jurka},
{\it Discovering simple DNA sequences by the algorithmic significance method},
Comput. Appl. Biosci., 9, (1993), p 407--411. 


\bibitem[Mon]{Mon}  T. {\sc Monteil},
{\it Estimating the number of ergodic measures of a minimal subshift knowing the
geometry of its Rauzy graphs},
preprint.




\end{thebibliography}
\end{document}